\documentclass[11pt]{amsart}
\usepackage{a4wide}
\usepackage{amsmath}
\usepackage[utf8]{inputenc}
\usepackage{amssymb}
\usepackage{bbm}
\usepackage{amsopn}
\usepackage{epsfig}
\usepackage{amsfonts}
\usepackage{latexsym}
\usepackage{graphicx}
\usepackage{enumerate}
\usepackage{color}
\usepackage{tikz}
\usepackage[colorlinks,linkcolor=blue,anchorcolor=blue,citecolor=blue]{hyperref}
\usepackage{cite}
\newtheorem{theorem}{Theorem}[section]
\newtheorem{lemma}[theorem]{Lemma}

\theoremstyle{definition}

\newtheorem{conjecture}[theorem]{Conjecture}

\theoremstyle{remark}

\numberwithin{equation}{section}

\newcommand{\R}{\ensuremath{\mathbb{R}}}
\newcommand{\N}{\ensuremath{\mathbb{N}}}
\newcommand{\Z}{\ensuremath{\mathbb{Z}}}
\newcommand{\Q}{\ensuremath{\mathbb{Q}}}

\newcommand{\A}{\ensuremath{\mathcal{A}}}
\newcommand{\B}{\ensuremath{\mathcal{B}}}

\newcommand{\set}[1]{\left\{#1\right\}}

\newcommand{\ep}{\varepsilon}
\newcommand{\f}{\infty}

\newcommand{\al}{\alpha}

\newcommand{\ggcd}{\mathbf{gcd}}

\begin{document}

\title{Rational points in translations of the Cantor set}

\author{Kan Jiang}
\address[K. Jiang]{Department of Mathematics, Ningbo University, Ningbo, People's Republic of China}
\email{jiangkan@nbu.edu.cn}

\author{Derong Kong}
\address[D. Kong]{College of Mathematics and Statistics, Center of Mathematics, Chongqing University, Chongqing 401331, People's Republic of China}
\email{derongkong@126.com}

\author{Wenxia Li}
\address[W. Li]{School of Mathematical Sciences, Key Laboratory of MEA (Ministry of Education) \& Shanghai Key Laboratory of PMMP, East China Normal University, Shanghai 200241, People's Republic of China}
\email{wxli@math.ecnu.edu.cn}

\author{Zhiqiang Wang}
\address[Z. Wang]{School of Mathematical Sciences, Key Laboratory of MEA (Ministry of Education) \& Shanghai Key Laboratory of PMMP, Shanghai 200241, People's Republic of China}
\email{zhiqiangwzy@163.com}

\date{\today}

\begin{abstract}
Given two coprime integers $p\ge 2$ and $q \ge 3$, let ${D_p\subset[0,1)}$ consist of all rational numbers which have a finite $p$-ary expansion, and let \[
 K(q, \A)=\set{\sum_{i=1}^\f\frac{d_i}{q^i}: d_i\in\A~ \forall i\in\N},
 \]
 where $\A\subset\set{0,1,\ldots, q-1}$ with cardinality $1<\#\A< q$.
In 2021 Schleischitz \cite{Schleischitz-2021} showed that $\#(D_p\cap K(q,\A))<+\f$. In this paper we show that for any $r\in\Q$ and for any $\al\in\R$,
\[
\#\big((r D_p+\al)\cap K(q,\A)\big)<+\f.
\]
\end{abstract}

\maketitle

\section{Introduction}

Given $q\in\N_{\ge 3}$ and $\A\subset\set{0,1,\ldots, q-1}$ with cardinality $1<\#\A<q$, we define the Cantor set {$K(q, \A)\subset [0,1]$} by
\[ K(q, \A):=\set{\sum_{i=1}^\f\frac{d_i}{q^i}: d_i\in\A~\forall {i\in\N}}. \]
In particular, the choice of $q=3$ and $\A=\{0,2\}$ corresponds to the classical middle-third Cantor set.
Given $p \in \N_{\ge 2}$, we define
\begin{equation}\label{def:D-p}
  D_p := \set{ \sum_{i=1}^{n} \frac{d_i}{p^i}:  d_i \in \set{0,1,\ldots, p-1} ~\forall 1\le i \le n;~ n \in \N}.
\end{equation}
Then $D_p$ consists of all rational numbers in ${[0,1)}$ which have a finite $p$-ary expansion. It is clear that $D_p$ is countably infinite and dense in $[0,1]$.

When $p=10, q=3$ and $\A=\set{0,2}$,
Wall \cite{Wall-1990} showed that
$$D_{10} \cap K(3,\set{0,2}) = \set{ \frac{1}{4}, \frac{3}{4}, \frac{1}{10}, \frac{3}{10}, \frac{7}{10}, \frac{9}{10}, \frac{1}{40}, \frac{3}{40}, \frac{9}{40}, \frac{13}{40}, \frac{27}{40}, \frac{31}{40}, \frac{37}{40}, \frac{39}{40}  }.$$
Later, Nagy \cite{Nagy-2001} proved that for each prime number $p \in\N_{\ge 4}$, the set $D_p \cap K(3,\set{0,2})$ is finite.
Bloshchitsyn \cite{Bloshchitsyn-2015} generalized this result and proved that if $p > q^2$ is a prime number, then the set $D_p \cap K(q,\A)$ is finite.
The general result was recently obtained by Schleischitz \cite[Corollary 4.4]{Schleischitz-2021} (some further extensions can be found in \cite{Shparlinski-2021} and \cite{Li-Li-Wu-2022}).

\begin{theorem}[Schleischitz, 2021]\label{thm:Schleischitz-2021}
  Let $p \in \N_{\ge 2}$ and $q\in\N_{\ge 3}$ with $\ggcd(p,q) =1$. If $\A\subset\set{0,1,\ldots, q-1}$ with $1<\#\A<q$,  then we have $$\# \big( D_p \cap K(q,\mathcal{A}) \big)<+\f.$$
\end{theorem}

In this paper we extend Theorem \ref{thm:Schleischitz-2021} as follows.
\begin{theorem}\label{thm:translation}
 Let $p \in \N_{\ge 2}$ and $q\in\N_{\ge 3}$ with $\ggcd(p,q) =1$. If $\A\subset\set{0,1,\ldots, q-1}$ with $1<\#\A<q$,  then for any $r \in \Q$ and any $\al\in\R$ we have
 $$\# \big( (r D_p +\al) \cap K(q,\mathcal{A}) \big)<+\f.$$
\end{theorem}

Note that in Theorem \ref{thm:Schleischitz-2021} the intersection $D_p\cap K(q,\A)$ contains only rational numbers, while in Theorem \ref{thm:translation} the intersection $(rD_p+\al)\cap K(q, \A)$ involves irrational numbers if $\al\notin\Q$.
To prove Theorem \ref{thm:translation}, we may assume that $\#\A=q-1$, which means the set $\A$ only misses one digit in $\{0,1,\ldots, q-1\}$. 
In terms of Theorem \ref{thm:translation} we make the following conjecture, which claims that the conclusion still holds also for irrational $x$.
\begin{conjecture}
Under the same condition as in Theorem \ref{thm:translation}, the conclusion
\[ \# \big( (rD_p + \alpha)\cap K(q,\A) \big) < +\f \]
also holds for any $r \notin \Q$ and any $\al\in\R$. 
\end{conjecture}
We will prove our main result Theorem \ref{thm:translation} in the next section.

\section{Proof of Theorem \ref{thm:translation}}\label{sec:translation}

In the following we fix two coprimes $p\in\N_{\ge 2}, q\in\N_{\ge 3}$ and the digit set $\A\subset\set{0,1,\ldots, q-1}$ with $\#A=q-1$.
For a real number $x \in \R$, we write $\langle x \rangle$ for the fractional part of $x$, i.e., $\langle x \rangle \in [0,1)$ and $x - \langle x \rangle \in \Z$.
For $x \in [0,1)$, the $q$-ary expansion of $x$ is the sequence $(x_i)$ in $\set{0,1,\ldots, q-1}^\N$ such that
\[
x=\sum_{i=1}^\f\frac{x_i}{q^i}.
\]
The $q$-ary expansion is unique except for countably many points that have precisely two $q$-ary expansions, one is finite and the other one ends in a periodic sequence with period $q-1$.
For convenience, for these countably many exceptional points, the $q$-ary expansion refers to the finite expansion. 
In order to prove Theorem \ref{thm:translation} we need the following lemma which can be deduced from Theorem \ref{thm:Schleischitz-2021}.

\begin{lemma}
  \label{lem:outside K(q,A)}
  Let $d_1 d_2\ldots d_k\in\{0,1,\ldots, q-1\}^k$ be a block and $r\in \Q \setminus \{0\}$. Then for all but finitely many $x \in r D_p$, the block $d_1 d_2\ldots d_k$ occurs in the $q$-ary expansion of $\langle x \rangle$ infinitely often.
\end{lemma}
\begin{proof}
  Write $r = s/t$ with $s \in \Z$, $t \in \N$, and $\ggcd(s,t) =1$.
  We can find $\ell \in \N$ such that
  \begin{equation}\label{def-ell}
    \ggcd\Big( \frac{t}{\ggcd(t,q^\ell)}, q \Big) =1.
  \end{equation}
  Let $h: = t / \ggcd(t,q^\ell)$.
 Since $\ggcd(h, q)=\ggcd(p,q)=1$, we have $\ggcd(h p, q^k) =1$. By Theorem \ref{thm:Schleischitz-2021}, we have
  \begin{equation}\label{eq:1}
  \# \big( D_{h p}\cap K(q^k, \B) \big) < +\infty,
  \end{equation}
  where
  \[
  \B = \big\{ 0,1,\ldots, q^k-1 \big\} \setminus \big\{ d_1 q^{k-1}+d_2 q^{k-2}+\cdots +d_{k-1}q+d_k \big\}.
  \]
  This implies that for any $y \in D_{hp} \setminus K(q^k,\B)$, the block $d_1 d_2\ldots d_k$ occurs in the $q$-ary expansion of $y$.
  Note that each $y \in D_{hp}$ has a purely periodic $q$-ary expansion because $\ggcd(h p, q) =1$ (cf. \cite[Proposition 2.1.2]{Dajani-Kraaikamp-2002}).
  Thus, for any $y \in D_{hp} \setminus K(q^k,\B)$, the block $d_1 d_2\ldots d_k$ occurs in the $q$-ary expansion of $y$ infinitely often.

  For $\ell$ defined in \eqref{def-ell}, consider the function $f$ defined by
  \[ f:\quad \R\to[0,1);\quad x\mapsto \langle q^\ell x \rangle. \]
  Note that the $q$-ary expansions of $\langle x \rangle$ and $f(x)$ have the same tail.
  Then for each $x \in f^{-1}\big( D_{hp} \setminus K(q^k,\B) \big)$, the block $d_1 d_2\ldots d_k$ occurs in the $q$-ary expansion of $\langle x \rangle$ infinitely often.
  It suffices to show that
  \begin{equation}\label{eq:finiteness}
    \# \Big( (rD_p) \setminus f^{-1}\big( D_{hp} \setminus K(q^k,\B) \big) \Big) < + \f.
  \end{equation}

  By (\ref{def:D-p}), we can rewrite $D_p$ as \[ D_p = \bigcup_{n=1}^\f \Big\{ \frac{d}{p^n}: d \in \{0,1,\cdots, p^n-1\}  \Big\}. \]
  Note that $h= t / \ggcd(t,q^\ell)$. Then for any $d/p^n \in D_p$, we have
  \begin{align*}
    f\Big( r \cdot \frac{d}{p^n} \Big)
    & = \Big\langle q^\ell \cdot \frac{s}{t} \cdot \frac{d}{p^n} \Big\rangle \\
    & = \Big\langle \frac{q^\ell}{\ggcd(t,q^\ell)} \cdot \frac{s}{h} \cdot \frac{d}{p^n} \Big\rangle \\
    & = \Big\langle \frac{q^\ell}{\ggcd(t,q^\ell)}  \cdot \frac{s h^{n-1} d}{(hp)^n} \Big\rangle \in D_{hp}.
  \end{align*}
  So, we obtain that $f(rD_p) \subset D_{hp}$, and then $r D_p\subset f^{-1}(D_{hp})$. 
 This implies that
 \begin{align*}
    (rD_p) \setminus f^{-1}\big( D_{hp} \setminus K(q^k,\B) \big) &= \Big((rD_p) \cap f^{-1}(D_{hp})\Big) \setminus f^{-1}\big( D_{hp} \setminus K(q^k,\B) \big)\\
   &= (rD_p) \cap \Big( f^{-1}(D_{hp}) \setminus f^{-1}\big( D_{hp} \setminus K(q^k,\B)\big) \Big)\\
   &= (rD_p) \cap f^{-1}\big( D_{hp}\cap K(q^k,\B) \big).
 \end{align*}
 Next, we show that $f$ is finite-to-one on $rD_p$, i.e., for any $y \in D_{hp}$, $\#\big( (rD_p) \cap f^{-1}(\{y\}) \big) < + \f$. Then, by (\ref{eq:1}) we obtain (\ref{eq:finiteness}) as desired.

  Let $p_1, p_2, \ldots, p_u$ be all distinct prime factors of $p$. Then each point in $D_p \setminus \{0\}$ has the unique representation of the form $$\frac{a}{p_1^{n_1} p_2^{n_2} \cdots p_u^{n_u}}$$ with $a \in \N$, $\ggcd(a,p) = 1$ and $n_1, n_2, \ldots, n_u$ are all nonnegative integers.
  For any $y = d/(hp)^n \in D_{hp}$, if
  \begin{align*}
    f\Big( r \cdot \frac{a}{p_1^{n_1} p_2^{n_2} \cdots p_u^{n_u}} \Big)
    & = \Big\langle q^\ell \cdot \frac{s}{t} \cdot \frac{a}{p_1^{n_1} p_2^{n_2} \cdots p_u^{n_u}} \Big\rangle \\
    & = \Big\langle \frac{q^\ell}{\ggcd(t,q^\ell)} \cdot \frac{s a}{h p_1^{n_1} p_2^{n_2} \cdots p_u^{n_u}} \Big\rangle \\
    & = y = \frac{d}{(hp)^n},
  \end{align*}
   then by using  $\ggcd(a,p)=1$ and $\ggcd(q,p)=1$ we obtain that
  $$p_1^{n_1} p_2^{n_2} \cdots p_u^{n_u} \mid s(hp)^n.$$
  This implies that all possible $n_1, n_2, \ldots, n_u$ are bounded.
  Thus, we conclude that $\#\big( (rD_p) \cap f^{-1}(\{y\}) \big) < + \f$,
  completing the proof. 
\end{proof}

We also need the following lemma.
\begin{lemma}\label{lem:induction}
  Let $\alpha \in \R$ and $d\in\set{1,2,\ldots, q-2}$ so that $q-d \not \in \A$. If the $q$-ary expansion $(\al_i)$ of $\langle \alpha \rangle$ satisfies $\al_i\ge d$ for all $i \in \N$, then
  \begin{equation}\label{eq:ineq-1}
    \#\big( (r D_p+\al)\cap K(q, \A) \big)<+\f
  \end{equation}for any $r \in \Q$.
\end{lemma}

\begin{proof}
If $r=0$, then (\ref{eq:ineq-1}) holds trivially. In the following we fix $r\in \Q \setminus \{0\}$.
Take $\al\in\R$ such that the $q$-ary expansion $(\al_i)$ of $\langle \alpha \rangle$ satisfies $\al_i\in\set{d, d+1, \ldots, q-1}$ for all $i\in\N$. We will prove (\ref{eq:ineq-1}) according to different properties of $(\al_i)$.
More precisely, for $0\le k\le q-d$ let
\begin{align*}
A_k:=&\Big\{ (a_i)\in\set{d, d+1,\ldots, q-1}^\N: \textrm{the length of any block in }(a_i)\textrm{ with }\\
& \hspace{4.2cm}\textrm{ each digit } < q-k \textrm{ is uniformly bounded} \Big\}.
\end{align*}
Then $A_0=\emptyset$ and $A_{q-d}=\set{d, d+1,\ldots, q-1}^\N$.
Note that $A_k \subset A_{k+1}$ for $0 \le k < q-d$. Thus,
$$(A_{q-d} \setminus A_{q-d-1}) \cup \cdots \cup (A_2 \setminus A_1) \cup (A_1 \setminus A_0) =  A_{q-d} \setminus A_0 = \{ d, d+1, \ldots, q-1\}^\N.$$
So it suffices to prove that for any $0\le k<q-d$,
\begin{equation}
  \label{eq:ineq-2}
  \#\big( (r D_p+\al)\cap K(q, \A) \big)<+\f\quad\textrm{if}\quad (\al_i)\in A_{k+1} \setminus A_{k}.
\end{equation}
We will split the proof of (\ref{eq:ineq-2}) into two cases: $0\le k<d$, and $d\le k<q-d$ (assuming $2d<q$).

\vspace{0.5em}
\noindent
\textbf{Case 1}. $(\al_i)\in {A_{k+1} \setminus A_k}$ for some $0\le k < d$.
Since $(\al_i)\in A_{k+1}$, the set
\begin{equation}\label{eq:ni}
\{ n_1,n_2, n_3,\ldots \} := \{ i \in \N: \alpha_i \ge q-k-1\}
\end{equation}
is infinite, where $n_1<n_2<n_3<\cdots$, and there exists $m\in\N$ such that $n_{i+1} - n_i \le m$ for all $i\in\N$.
Note that $0\le k < d < q$. Then $1\le q-d+k \le q-1$.

Take $x \in r D_p$ such that the block $(q-d+k)^{m+1}0$ occurs in the $q$-ary expansion $(x_i)$ of $\langle x \rangle$ infinitely often,
and let
\begin{equation}\label{eq:ki}
\{ k_1, k_2,k_3,\ldots \} := \{ i \in \N_{\ge n_1}: x_{i+1} x_{i+2} \cdots x_{i+m+2} = (q-d+k)^{m+1}0 \},
\end{equation}
where $k_1<k_2<k_3<\cdots$.
Note that the $q$-ary expansion of $\langle x \rangle$ is eventually periodic. Then the sequence $\{ k_{i+1}-k_i\}_{i \in \N}$ is bounded.

\noindent{\bf Claim}: there exists $i_0 \in \N$ such that $\alpha_{k_{i_0} +j} \le q-k-1$ for all $1 \le j \le m+ 2$.

\noindent
Suppose on the contrary that {the claim} fails. Then for each $i \in \N$  the block $\alpha_{k_i+1} \alpha_{k_i+2}  \cdots \alpha_{k_i+m+2}$ contains at least one digit $\ge q-k$.
Note that $\{ k_{i+1}-k_i\}_{i \in \N}$ is bounded. Then the length of any block in $(\al_i)$ with each digit $<q-k$ is uniformly bounded. This implies that $(\al_i)\in A_k$, leading to a contradiction with our assumption.

Since the sequence $\{ n_{i+1} - n_i \}_{i \in \N}$ is bounded by $m$ and $k_{i_0} \ge n_1$, there exists $1 \le m' \le m$ such that $k_{i_0}+m'\in\set{n_i}$, and thus by (\ref{eq:ni}) and {the claim} it follows that $\alpha_{k_{i_0} + m'} = q-k-1.$
Thus, by (\ref{eq:ki}) and using $\alpha_i \ge d$ for all $i \in \N$ we obtain that
\begin{align*}
 \sum_{i=1}^{\f} \frac{x_{k_{i_0} +m' -1 +i}}{q^i} + \sum_{i=1}^{\f} \frac{\alpha_{k_{i_0} +m' -1 + i}}{q^i} &>\left( \frac{q-d+k}{q} + \frac{q-d+k}{q^2}\right) + \left(\frac{q-k-1}{q} + \frac{d}{q^2}\right) \\
 &\ge 1+ \frac{q-d}{q},
 \end{align*}
 and
 \begin{align*}
  \sum_{i=1}^{\f} \frac{x_{k_{i_0} +m' -1 +i}}{q^i} + \sum_{i=1}^{\f} \frac{\alpha_{k_{i_0} +m' -1 + i}}{q^i} &<\left(\frac{q-d+k}{q} + \frac{q-d+k+1}{q^2}\right) + \frac{q-k}{q} \\
  &\le 1 + \frac{q+1-d}{q}.
  \end{align*}
Whence,
\begin{align*}
\big\langle  q^{k_{i_0}+m'-1} (x+\alpha) \big\rangle &=\Big\langle \sum_{i=1}^{\f} \frac{x_{k_{i_0} +m' -1 +i}}{q^i} + \sum_{i=1}^{\f} \frac{\alpha_{k_{i_0} +m' -1 + i}}{q^i}\Big\rangle \\
&\in \Big( \frac{q-d}{q}, \frac{q+1-d}{q} \Big).
\end{align*}
This implies that $x+ \alpha \not\in K(q,\A)$.
So, applying Lemma \ref{lem:outside K(q,A)} for the block $(q-d+k)^{m+1}0$ involved in the $q$-ary expansion of $x$ in (\ref{eq:ki}), it follows that $ \#\big( (rD_p +\alpha) \cap K(q,\A) \big) < +\f$.

\vspace{0.5em}
\noindent
\textbf{Case 2}. $(\al_i)\in {A_{k+1} \setminus A_k}$ for some $d\le k < q-d$ with the assumption $2d<q$.

Similar to Case 1, since $(\al_i)\in A_{k+1}$, the set
$$ \{ n_1, n_2 ,n_3, \ldots \}: = \{ i \in \N: \alpha_i \ge q-k-1\} $$
is infinite, where $n_1< n_2 <n_3<\cdots$, and there exists $m\in\N$ such that $n_{i+1} - n_i \le m$ for all $i\in\N$.
Note that $1\le d\le k<q-d$. Then $1\le k-d+1\le q-1$.
Take $x\in r D_p$ such that the block $(k-d+1)^{m}0$ occurs in the $q$-ary expansion $(x_i)$ of $\langle x \rangle$ infinitely often, and write
\begin{equation}\label{eq:ki-Case2}
  \{ k_1 , k_2 , \ldots \} := \{ i \in \N_{\ge n_1}: x_{i+1} x_{i+2} \cdots x_{i+m+1} = (k-d+1)^{m}0 \},
\end{equation}
where $k_1 < k_2 < \ldots$.
Note that the $q$-ary expansion of $\langle x \rangle$ is {eventually periodic}. Then the sequence $\{ k_{i+1}-k_i\}_{i \in \N}$ is bounded. Since $(\al_i)\notin A_k$, by the same argument as in Case 1 we can find $i_0\in\N$ such that $\alpha_{k_{i_0} +j} \le q-k-1$ for all $1 \le j \le m+ 1$, and there exists $1 \le m' \le m$ such that $k_{i_0} + m' \in \{n_i\}$ and $\alpha_{k_{i_0} + m'} = q-k-1$. Thus,
\[ \sum_{i=1}^{\f} \frac{x_{k_{i_0} +m' -1 +i}}{q^i} + \sum_{i=1}^{\f} \frac{\alpha_{k_{i_0} +m' -1 + i}}{q^i} > \frac{k-d+1}{q} + \frac{q-k-1}{q} = \frac{q-d}{q}, \]
and
\begin{align*}
  &\; \sum_{i=1}^{\f} \frac{x_{k_{i_0} +m' -1 +i}}{q^i} + \sum_{i=1}^{\f} \frac{\alpha_{k_{i_0} +m' -1 + i}}{q^i} \\
  < &\; \Bigg( \sum_{i=1}^{m-m'+1} \frac{k-d+1}{q^i} + \frac{1}{q^{m-m'+2}} \Bigg) + \Bigg( \sum_{i=1}^{m-m'+2} \frac{q-k-1}{q^i} + \frac{1}{q^{m-m'+2}} \Bigg)\\
  = &\; \sum_{i=1}^{m-m'+1} \frac{q-d}{q^i} + \frac{q-k+1}{q^{m-m'+2}} \\
  \le &\; \frac{q+1 - d}{q}.
\end{align*}
Whence, \[ \big \langle q^{k_{i_0} +m'-1} (x+\alpha) \big \rangle \in \Big( \frac{q-d}{q}, \frac{q+1-d}{q} \Big), \]
which implies $ x+ \alpha \not\in K(q,\A)$. Therefore, applying Lemma \ref{lem:outside K(q,A)} for the block $(k-d+1)^{m}0$ involved in the $q$-ary expansion of $x$ in (\ref{eq:ki-Case2}), we deduce  that $ \#\big( (rD_p +\alpha) \cap K(q,\A) \big) < +\f$. This together with Case 1 proves (\ref{eq:ineq-2}), completing the proof.
\end{proof}

Now we are ready to prove Theorem \ref{thm:translation}.
\begin{proof}[Proof of Theorem \ref{thm:translation}]
  If $r=0$, then it holds trivially. In the following we fix $r \in \Q \setminus \{0\}$ and $\alpha \in \R$. Let $(\alpha_i)$ denote the $q$-ary expansion of $\langle\alpha\rangle$.
  Since $\#\A = q-1$, there exists $1 \le d \le q$ such that $q-d \not\in \A$.
  We split the proof  into two cases according to $(\al_i)$.

  \noindent\textbf{Case 1}. For any $m \in \N$ the sequence $(\alpha_i)$ contains the block $0^m$.

  Fix $x \in rD_p$ such that the block $(q-d)0$ occurs in the $q$-ary expansion $(x_i)$ of $\langle x \rangle$ infinitely often, and write
  \begin{equation}\label{eq:index-ki}
    \{ k_1 , k_2 , \ldots \} := \{ i \in \N: x_{i+1} x_{i+2} =(q-d)0 \},
  \end{equation}
  where $k_1 < k_2 < \ldots$.
  Note that the $q$-ary expansion of $\langle x \rangle$ is eventually periodic.
  Then the sequence $\{ k_{i+1}-k_i\}_{i \in \N}$ is bounded. Observe that $(\al_i)$ contains arbitrarily long consecutive zeros. There   exists $i_0 \in \N$ such that $\alpha_{k_{i_0}+1} \alpha_{k_{i_0}+2} = 00$.
  Thus, by using
 $x_{k_{i_0}+1} x_{k_{i_0}+2}= (q-d)0$  it follows that
  \[ \frac{q-d}{q} < \sum_{i=1}^{\f} \frac{x_{k_{i_0} +i}}{q^i} + \sum_{i=1}^{\f} \frac{\alpha_{k_{i_0} + i}}{q^i} < \frac{q+1-d}{q}, \]
 which implies $\big\langle q^{k_{i_0}} (x+\alpha) \big\rangle \not\in K(q,\A).$
  So, $ x+ \alpha \not\in K(q,\A)$. Applying Lemma \ref{lem:outside K(q,A)} for the block $(q-d)0$ involved in the $q$-ary expansion of $x$ in (\ref{eq:index-ki}), it follows that $\#\big( (rD_p +\alpha) \cap K(q,\A) \big) < +\f$.

  \noindent\textbf{Case 2}. There exists $m \in \N$ such that the sequence $(\alpha_i)$ does not contain the block $0^m$.

  Note that  $K(q,\A) = K(q',\A')$, where \[q' = q^{m+1}\quad\textrm{and}\quad \A' = \set{ \sum_{i=0}^m \ep_i q^i: \ep_i \in \A  ~\forall 0\le i\le m }. \]
  Then it suffices to show that $\# \big( (r D_p+\alpha) \cap K(q',\A') \big)<+\f$.
  Since $\ggcd(p,q)=1$, we have $\ggcd(p,q')=1$.
  Observe that \[ \langle \alpha \rangle = \sum_{i=1}^{\f} \frac{\alpha_i}{q^i} = \sum_{k=1}^{\f} \frac{1}{q^{k(m+1)}}\left( \sum_{i=0}^{m} \alpha_{k(m+1)-i} q^i\right) = \sum_{k=1}^{\f} \frac{\alpha_k'}{(q')^k}, \]
  where
  \begin{equation}\label{def:alpha-k-prime}
    \alpha_k' := \sum_{i=0}^{m} \alpha_{k(m+1)-i} q^i.
  \end{equation}
  Then the sequence $(\alpha_k')$ is the $q'$-ary expansion of $\langle\alpha\rangle$.
  Since the sequence $(\alpha_i)$ does not contain the block $0^m$, for each $k\in \N$ there exists $1 \le i_k \le m$ such that $\alpha_{k(m+1)-i_k} >0$.
  It follows from (\ref{def:alpha-k-prime}) that $\alpha_k' \ge q \ge d$ for all $k \in \N$.
  Furthermore, since $q-d \not\in \A$, we have \[ q' - d = q-d + (q-1) q + (q-1) q^2 + \cdots + (q-1) q^{m}\not\in \A'.\]
  Clearly, we have $1\le d \le q < q'-1$.
  Thus, applying Lemma \ref{lem:induction} for the set $K(q',\A')$ we conclude that $\# \big( (rD_p+\alpha) \cap K(q',\A') \big)<+\f$.
\end{proof}

\section*{Acknowledgements}
The authors thank the anonymous referees for many useful suggestions which improve the
presentation of the paper. 
K.~Jiang was supported by NSFC No.~11701302,  Zhejiang Provincial NSF No.~LY20A010009, and the K.C. Wong Magna Fund in Ningbo University. D.~Kong was supported by NSFC No.~11971079. W.~Li was  supported by NSFC No.~12071148 and Science and Technology Commission of Shanghai Municipality (STCSM) No. 22DZ2229014.  
Z. Wang was supported by the Fundamental Research Funds for the Central Universities No.~YBNLTS2023-016.


\begin{thebibliography}{10}

\bibitem{Bloshchitsyn-2015}
V. Ya. Bloshchitsyn.
\newblock Rational points in m-adic Cantor sets.
\newblock {\em J. Math. Sci. (N.Y.)}, 211(6):747–751, 2015.


\bibitem{Dajani-Kraaikamp-2002}
K. Dajani, C. Kraaikamp,
\newblock \emph{Ergodic theory of numbers}.
\newblock Carus Mathematical Monographs, 29. Mathematical Association of America, Washington, DC, 2002.


\bibitem{Li-Li-Wu-2022}
B.~Li, R.~Li and Y.~Wu.
\newblock Rational numbers in $\times b$-invariant sets.
\newblock {\em Proc. Amer. Math. Soc.}, 151(5), 1877–1887, 2023.


\bibitem{Nagy-2001}
J. Nagy.
\newblock Rational points in Cantor sets.
\newblock {\em Fibonacci Quart.}, 39(3):238–241, 2001.


\bibitem{Schleischitz-2021}
J. Schleischitz.
\newblock On intrinsic and extrinsic rational approximation to Cantor sets.
\newblock {\em Ergodic Theory Dynam. Systems}, 41(5):1560–1589, 2021.


\bibitem{Shparlinski-2021}
I. E. Shparlinski.
\newblock On the arithmetic structure of rational numbers in the Cantor set,
\newblock {\em Bull. Aust. Math. Soc.}, 103(1), 22–27, 2021.

\bibitem{Wall-1990}
C. R. Wall.
\newblock Terminating decimals in the Cantor ternary set.
\newblock {\em Fibonacci Quart.}, 28(2):98–101, 1990.

\end{thebibliography}
\end{document}